# A note on the modes of the negative binomial distribution of order *k*, type I

Costas Georghiou[a], Andreas N. Philippou[b] and Zaharias M. Psillakis[c,*]

[a] School of Engineering, University of Patras, 26500, Greece; [b] Department of Mathematics, University of Patras, 26500, Greece; [c] Department of Physics, University of Patras, 26500, Greece; [*] Corresponding author

c.georghiou@upatras.gr, anphilip@math.upatras.gr, psillaki@physics.upatras.gr

**ABSTRACT**

Upper and lower bounds are derived for the mode(s) $m_k(r,p)$ of the negative binomial distribution of order *k*, type *I*, with parameters *r* and *p*, say $NB_{k,I}(r,p)$, which are employed to establish an explicit formula for $m_k(r,p)$ in terms of *r* and *k* when *p = 1/2*. It is also shown as a direct consequence of the upper bound alone that $m_k(r,p) = k$ when *r = 1*. The derivation of the bounds is based on a known recurrence relation satisfied by the probability mass function of $NB_{k,I}(r,p)$.

1. **Introduction and summary**

The discrete distributions of order k have received a lot of attention during the last several decades [1 - 17], however, only partial results are known for the modes of the Poisson distribution of order k and the negative binomial distribution of the same order, type I. See Luo (1987), Georghiou, Philippou and Saghafi (2013), Philippou (2014), and Shao and Fu (2016). We reconsider presently the latter.

For any given positive integers *k* and *r* and $0 < p < 1$ ($q = 1 - p$), a random variable X is said to be distributed as negative binomial distribution of order k, type I, with parameters r and p, denoted $NB_{k,I}(r,p)$, if

(1.1) $\quad P_n = P(X = n) = p^n \sum \binom{n_1 + n_2 + \cdots + n_k + r - 1}{n_1, n_2, \ldots, n_k, r - 1} \left(\frac{q}{p}\right)^{n_1+n_2+\cdots+n_k},$

for $n \geq kr$ and *0* otherwise, where the summation is taken over all *k*-tuples of non-negative integers $n_1, n_2, \ldots, n_k$ such that $n_1 + 2n_2 + \cdots + kn_k = n - kr$ [2, 3, 12, 14, 16].

For $k = 1$ $NB_{k,I}(r,p)$ reduces to the negative binomial distribution with parameters $r$ and $p$, and for $r = 1$ it reduces to the geometric distribution of order k with parameter $p$.

We denote presently by $m \coloneqq m_k(r,p)$ the mode(s) of $NB_{k,I}(r,p)$, i.e. the value(s) of *n* for which $P_n$ attains its maximum. It is well-known that $m_1(r,p) = \lfloor (r-1)/p \rfloor + 1$, if $(r-1)/p$ is not an integer, and $m_1(r,p) = (r-1)/p$ and $(r-1)/p + 1$, if $(r-1)/p$ is an integer. Here and in the sequel, $\lfloor x \rfloor$ denotes the greatest integer not exceeding x.

In a recent paper, Shao and Fu (2016) obtained the following modes for special cases of $NB_{k,I}(r,p)$

(1.2a)  $m_k(1, p) = k$ $(k \geq 1, 0 < p < 1)$, and

(1.2b)  $m_2(2, 1/2) = 6,7,8$, $m_3(2, 1/2) = 16$, $m_2(3, 1/2) = 13$.

See their Theorems 2.6, 3.6, 3.7, and 3.9, as well as an application of the modes of $NB_{k,I}(r,p)$ in continuous sampling plans [see their Remark 3.11].

Using the following recurrence relation of Philippou and Georghiou (1989)

(1.3)  $(n-kr)P_n = q\sum_{j=1}^{k}(n-kr+j(r-1))p^{j-1}P_{n-j}$, $(n > kr, k, r \geq 1, 0 < p < 1)$,

with $P_n = 0$ for $n < kr$ and $P_{kr} = p^{kr}$ [see their Theorem 3.1], we show in Section 2 that

(1.4)  $m_k(r,p) \leq kr + \left\lfloor \frac{(r-1)[1-p^k(1+kq)]}{qp^k} \right\rfloor$, $k, r \geq 1$, $0 < p < 1$

and

(1.5)  $m_k(r,p) \geq kr + \lfloor \rho \rfloor$, $k, r \geq 2$, $0 < p \leq (r-1)/r$,

where $\rho$ depends on the roots of a certain quadratic polynomial [see Theorem 2.2].

By the definition of the mode and the sharp upper bound (1.4) we readily get $m_k(1,p) = k$. The derivation of the remaining three modes in (1.2) follows from (1.4) and some, very few, simple calculations using (1.3). We also show, by means of (1.4) and (1.5), that

(1.6)  $m_k(r, 1/2) = kr + (r-1)(2^{k+1} - k - 2) - 1 + \delta_{k,r}$, $k, r \geq 2$,

where $\delta_{k,r} = 0, \pm 1$ if $k = r = 2$, and 0 otherwise.

In the sequel, unless otherwise explicitly stated, k and r are positive integers and $0 < p < 1$.

2. Main results

In the present section, we state and prove the following three theorems.

**Theorem 2.1:** *The mode* $m := m_k(r,p)$ *of the negative binomial distribution of order k, type I, with parameters r and p satisfies the inequality*

$m \leq kr + \left\lfloor \frac{(r-1)[1-p^k(1+kq)]}{qp^k} \right\rfloor := \overline{m}_k(r,p) := \overline{m}$. □

Taking $r = 1$ in the theorem we readily get

**Corollary 2.1:** *The mode m of the geometric distribution of order k with parameter p is $m = k$.* □

Furthermore, noting that $P_m = \max\{P_{kr}, \ldots, P_{\overline{m}}\}$, by the definition of the mode and Theorem 2.1, and using Theorem 2.1 and the effective recurrence relation (1.3) to calculate $\overline{m}$ and $P_{kr}, \ldots, P_{\overline{m}}$ for given $k, r$, and $p$, we can evaluate the mode $m$ for these values. In particular,

**Corollary 2.2:** (a) $m_2(2, 1/2) = 6, 7, 8$, (b) $m_3(2, 1/2) = 16$, and (c) $m_2(3, 1/2) = 13$. □

Furthermore, by means of a simple computer program, based on the definition of the modes, and a laptop we get Table 1. The modes $m_k(r, 0.5)$ for k, r = 2, 3, 4, 5 led to Theorem 2.3. The modes $m_k(r, p) = kr$ for k, r = 2, 3, 4, 5 and p = 0.9, 0.95, 0.99 indicate that $m_k(r, p) = kr$ for k, r = 2, 3, 4, 5 and p ≥ 0.9.

**Table 1:** Modes $m_k(r, p)$ of $NB_k(r, p)$ for k, r = 2, 3, 4, 5 and p = 0.5, 0.6, 0.7, 0.8, 0.9, 0.95, 0.99.

| $m_k$(r, 0.5) | | | | |
|---|---|---|---|---|
| k/r | 2 | 3 | 4 | 5 |
| 2 | 6, 7, 8 | 13 | 19 | 25 |
| 3 | 16 | 30 | 44 | 58 |
| 4 | 33 | 63 | 93 | 123 |
| 5 | 66 | 128 | 190 | 252 |
| $m_k$(r, 0.6) | | | | |
| k/r | 2 | 3 | 4 | 5 |
| 2 | 4 | 10 | 14 | 19 |
| 3 | 11 | 20 | 29 | 38 |
| 4 | 19 | 36 | 53 | 70 |
| 5 | 33 | 63 | 93 | 122 |
| $m_k$(r, 0.7) | | | | |
| k/r | 2 | 3 | 4 | 5 |
| 2 | 4 | 8 | 12 | 15 |
| 3 | 6 | 15 | 21 | 27 |
| 4 | 8 | 24 | 34 | 45 |
| 5 | 20 | 36 | 53 | 69 |
| $m_k$(r, 0.8) | | | | |
| k/r | 2 | 3 | 4 | 5 |
| 2 | 4 | 6 | 10 | 12 |
| 3 | 6 | 9 | 15 | 21 |
| 4 | 8 | 12 | 24 | 31 |
| 5 | 10 | 24 | 34 | 44 |
| $m_k(r, p) = kr$ for k, r = 2, 3, 4, 5 and p = 0.9, 0.95, 0.99. | | | | |

**Theorem 2.2:** The mode $m := m_k(r, p)$ of the negative binomial distribution of order k, type I, with parameters r and p satisfies the inequality

$$m \geq kr + \lfloor \rho \rfloor := \underline{m_k}(r, p) := \underline{m}, \quad k, r \geq 2, \quad 0 < p \leq (r-1)/r,$$

where $\rho = v_2$, the largest real root of the quadratic polynomial

$$f(v) = -qp^k v^2 + \left(p^k[q - (kq + q + 1)(r-1)] + r - 1\right)v +$$

$$+ (r-1)\left\{(r-1)\left(1 - p^k(1 + kq)\right) + p^k(k+q) - (1-p^k)/q\right\},$$

if $v_2 > k$, or $\rho = k$ if $v_2 \leq k$. □

**Theorem 2.3:** *The mode $m_k(r, 1/2)$ of the negative binomial distribution of order k, type I, with parameters r and p = ½ is given by the formula*

$$m_k(r, 1/2) = kr + (r-1)(2^{k+1} - k - 2) - 1 + \delta_{k,r}, \quad k, r \geq 2,$$

where $\delta_{k,r} = 0, \pm 1$ if $k = r = 2$, and 0 otherwise. □

**Proof of Theorem 2.1:** By the definition of the mode we have $P_n \leq P_m$ for any $n \geq kr$, which gives

$(n - kr)P_n \leq \left(q \sum_{j=1}^{k}(n - kr + j(r-1))p^{j-1}\right)P_m$, by (1.3),

$$= \left(q(n-kr)\sum_{j=1}^{k} p^{j-1} + q(r-1)\sum_{j=1}^{k} jp^{j-1}\right)P_m$$

$$= (n-kr)(1-p^k)P_m + (r-1)\frac{kp^{k+1} - (k+1)p^k + 1}{q}P_m$$

since $\sum_{j=1}^{k} p^{j-1} = (1-p^k)/(1-p)$ and $\sum_{j=1}^{k} jp^{j-1} = (kp^{k+1} - (k+1)p^k + 1)/(1-p)^2$. The theorem follows by taking $n = m$ and solving for $m$. □

**Remark 2.1:** It is well-known [14, 16] that if $q \to 0$ and $rq = \lambda_r \to \lambda\ (> 0)$ as $r \to \infty$, then the negative binomial distribution of order k, type I, with parameters r and p, properly shifted to the left, converges to the Poisson distribution of order k with parameter λ, say $P_k(\lambda)$, as $r \to \infty$. If $q \to 0$ and $rq = \lambda_r \to \lambda\ (> 0)$ as $r \to \infty$, then $\overline{m} - kr \to \left\lfloor \frac{1}{2}k(k+1)\lambda \right\rfloor$ as $r \to \infty$, which is the upper bound for the mode of $P_k(\lambda)$ derived by Georghiou, Philippou and Saghafi (2013). □

In order to prove Theorem 2.2, we need the following result.

**Lemma 2.1:** For $qr > 1$ and $kr \leq n \leq kr + k$ we have $P_n > P_{n-1}$.

**Proof of Lemma 2.1:** For $n = kr$ the lemma is true by the definition of $P_n$. For $n = kr + 1$ it is also true, since by the recurrence (1.3) and the assumption of the lemma $P_{kr+1} = qrP_{kr} > P_{kr}$. Finally, for $n = kr + v$, $1 < v \leq k$, we have from (1.3)

$$vP_{kr+v} = q[(v + r - 1)P_{kr+v-1} + (v + 2(r-1))pP_{kr+v-2} + \cdots + (v + v(r-1))p^{v-1}P_{kr}]$$
$$> q\{(v + r - 1)P_{kr+v-1}$$
$$+ p[(v - 1 + r - 1)P_{kr+v-2} + (v - 1 + 2(r-1))pP_{kr+v-3} + \cdots$$
$$+ (v - 1 + (v-1)(r-1))p^{v-2}P_{kr}]\} = q(v + r - 1)P_{kr+v-1} + p(v-1)P_{kr+v-1}$$
$$= (v + qr - 1)P_{kr+v-1} > vP_{kr+v-1,}$$

which completes the proof of the lemma. □

We note that if $qr = 1$ the lemma is still valid except for $n = kr + 1$ in which case $P_{kr+1} = P_{kr}$.

**Proof of Theorem 2.2:** Setting $v := n - kr$ and $\Delta_v = P_{kr+v} - P_{kr+v-1}$ ($v \geq 0$), we have $\Delta_0 = p^{kr} > 0$, $\Delta_1 = P_{kr+1} - P_{kr} = (qr - 1)p^{kr} > 0$, since $qr > 1$ by assumption (the case $qr = 1$ will be treated later on). Also, by Lemma 2.1, we have $\Delta_2 > 0, \Delta_3 > 0, \ldots, \Delta_k > 0$ meaning hat $P_{kr+v}$ is strictly increasing on the set $\{0, 1, \ldots, k\}$. The question then arises as to whether $P_{kr+v}$ keeps increasing (strictly) beyond $v = k$ and if so, up to which value $v = \rho$ one can say with certainty that $\Delta_{\lfloor \rho \rfloor} > 0$. This ρ will determine the lower bound for the mode. To do this we need some recursive relation on the $\Delta$'s. For this purpose, by subtracting $v\, P_{vkr+v-1}$ from both sides of (1.3), we get

(2.1) $$v\Delta_v = (-pv + q(r-1))P_{kr+v-1} + q\sum_{j=2}^{k}(v + j(r-1))p^{j-1}P_{kr+v-j}.$$

Now multiply the above equation by $(v-1)$ and replace $(v-1) P_{kr+v-1}$ with the help of (1.3) to get

$$v(v-1)\Delta_v = qp^{k-1}\left(-pv^2 + (p+(q-pk)(r-1))v + q(r-1)(kr-k-1)\right)P_{kr+v-1-k} +$$

$$+q(r-1)\sum_{j=1}^{k-1}(v-1+j(qr-1))p^{j-1}P_{kr+v-1-j}.$$

By setting, successively, $P_{kr+v-1-j} = \Delta_{v-1-j} + P_{kr+v-2-j}$, j = 1, 2, …, k-1, we finally get the recursion

(2.2) $$v(v-1)\Delta_v = f(v)P_{kr+v-1-k} + q(r-1)\sum_{j=1}^{k-1}(\sum_{i=1}^{j}(v-1+i(qr-1))p^{i-1})\Delta_{v-1-j},$$

where $f(v)$ is as in the theorem.

We immediately note that the coefficient of $\Delta_{v-1-j}$ in (2.2), is always positive and by Lemma 2.1 $\Delta_v > 0$ for $0 \leq v \leq k$. The positivity, then, of $\Delta_v$, for v>k, depends entirely on the properties of the quadratic polynomial $f(v)$. Since the coefficient of $v^2$ is negative, $f(v)$ takes positive values only between its real roots. But $f(1) = (qr-1)\frac{(r-1)[1-p^k(1+kq)]}{q}$ is positive, meaning that f(v) has indeed two real roots $v_1, v_2$ such that $v_1 < 1 < v_2$. It is easy then to see that $f(v) > 0$ for 0 <v<$v_2$ and $f(v) < 0$ for v>$v_2$ establishing, thus, the theorem for $0 < p < (r-1)/r$.

To complete the proof, we need to examine the case p = (r – 1)/r. Indeed, by setting q = 1/r in (2.2) we get, after removing the common factor $v-1$,

(2.3) $$v\Delta_v = \left((r-1)(1-p^{k-1}) - (\tfrac{v-1}{r}+kp)p^k\right)P_{kr+v-1-k} + (r-1)\sum_{j=1}^{k-1}(1-p^j)\Delta_{v-1-j},$$

where now the coefficient of $P_{kr+v-1-k}$ becomes negative when $v > r(r-1)\left(\frac{r}{r-1}\right)^k - k(r-1) - (r-1)(r+1)$ giving for $\rho$ the exact value $\rho = r(r-1)\left(\frac{r}{r-1}\right)^k - (r-1)(k+r+1)$. □

**Remark 2.2:** If $q \to 0$ and $rq = \lambda_r \to \lambda\ (>1)$ as $r \to \infty$, then $\rho = \underline{m} - kr \to \left\lfloor \tfrac{1}{2}k(k+1)(\lambda-1)\right\rfloor + 1$, as $r \to \infty$, which is the lower bound for the mode of $P_k(\lambda)$ derived by Georghiou, Philippou and Saghafi (2013). □

**Remark 2.3:** Equation (2.1) is valid for k = 1 giving $v\Delta_v = (-pv + q(r-1))P_{r+v-1}$. The coefficient of $P_{r+v-1}$ is negative if $v > q(r-1)/p$ and this gives $\underline{m}_1(r,p) = r + \lfloor q(r-1)/p \rfloor = 1 + \lfloor (r-1)/p \rfloor$. At the same time (1.4) gives the same value for $\overline{m}_1(r,p)$. However, if $v = q(r-1)/p$ is an integer then $\Delta_v = 0$ i.e. $P_{r+v} = P_{r+v-1}$ and, therefore $(r-1)/p$ is also a mode. □

**Proof of Theorem 2.3:** By assumption $k \geq 2, r \geq 2$ and $p = 1/2$. For $(k,r) = (2,2)$ the theorem holds true because of Corollary 2.2(a).

Let then $(k, r) \neq (2,2)$. By Theorems 2.1 and 2.2 and the assumptions we have

$$\underline{m} = kr + \rho = \leq m \leq \overline{m} = kr + (r-1)(2^{k+1} - k - 2).$$

We will show first that the root $\rho = \rho_k(r, 1/2)$ of $f(v)$ is such that $\lfloor \rho \rfloor = \overline{m} - kr - 1$, and hence

$$\overline{m} - 1 \leq m \leq \overline{m}.$$

Indeed, from the definition of $f(v)$ and the assumptions of Theorem 2.2 we have

$$f(\overline{m} - kr) = (r-1)[(k+3)/2^{k+1} - 1] < 0 \text{ and } f(\overline{m} - kr - 1) = 2^{-k}(r-2) \geq 0,$$

meaning that the root is either $\overline{m} - kr - 1$ (if r = 2) or lies between $\overline{m} - kr - 1$ and $\overline{m} - kr$ (if r > 2). In order to complete the proof, it suffices, then, to show that $P_{\overline{m}} < P_{\overline{m}-1}$. In fact, from the recurrence relation (1.3) we have, since $(\overline{m} - kr) > 0$ and $P_{\overline{m}-j} < P_{\overline{m}-1}$ $(2 \leq j \leq k)$,

$$(\overline{m} - kr)P_{\overline{m}} = q \sum_{j=1}^{k} ((\overline{m} - kr) + j(r-1))p^{j-1}P_{\overline{m}-j} < q \sum_{j=1}^{k} ((\overline{m} - kr) + j(r-1))p^{j-1}P_{\overline{m}-1}$$

$$= q \left( (\overline{m} - kr) \sum_{j=1}^{k} p^{j-1} + (r-1) \sum_{j=1}^{k} j p^{j-1} \right) P_{\overline{m}-1} = (\overline{m} - kr)P_{\overline{m}-1}.$$

The above inequality implies $P_{\overline{m}} < P_{\overline{m}-1}$, as claimed, and this completes the proof of the theorem. □

**Remark 2.4:** If k = r = 2 and p = 1 / 2, Theorems 2.1 and 2.2, respectively, give $\overline{m}_2(2, 1/2) = 4 + 4 = 8$ and $\underline{m}_2(2, 1/2) = 4 + 3 = 7$. Hence, $7 \leq m_2(2, 1/2) \leq 8$. But equation (2.3) gives $\Delta_3 = 0$, showing that $P_7 = P_6$, and $\Delta_4 = 0$, showing that $P_8 = P_7$. Therefore, $m_2(2, 1/2) = 6, 7, 8$. □

**Disclosure statement**

No potential conflict of interest was reported by the authors.